\documentclass[11pt]{article}

\usepackage{amssymb}
\usepackage{latexsym}
\usepackage{amsfonts}
\usepackage{amsmath}
\newtheorem{thrm}{Theorem}[section]
\newtheorem{lemma}[thrm]{Lemma}
\newtheorem{prop}[thrm]{Proposition}

\newtheorem{remark}[thrm]{Remark}
\newtheorem{example}[thrm]{Example}
\newtheorem{assumption}[thrm]{Assumption}
\numberwithin{equation}{section}

\usepackage[dvips]{color}

\def\wt{\widetilde}

\def\P{\mathbb{P} }
\def\E{\mathbb{E} }
\def\R{\mathbb{R} }
\def\V{\mathbb{V} }

  \def\sL {{\cal L}}

\def\R {{\mathbb R}}
\def\wt{\widetilde}
\def\wh{\widehat}

\def\E{{\mathbb E}}
\def\P{{\mathbb P}}

\def\bea{\begin{align*}}
\def\eea{\end{align*}}
\def\bee{\begin{equation}}
\def\eee{\end{equation}}

\def\eps{\varepsilon}

\def\wh{\widehat}

\begin{document}
\allowdisplaybreaks

\title{Strong law of large numbers for supercritical superprocesses under second moment condition}
\author{ \bf Zhen-Qing Chen\footnote{The research of ZQ is partially supported
by NSF grant DMS-1206276 and
 NSFC (Grant No. 11128101).}\hspace{1mm}\hspace{1mm}
 Yan-Xia
Ren\footnote{The research of YX is supported by NSFC (Grant No.  11271030 and 11128101).\hspace{1mm} } \hspace{1mm}\hspace{1mm}
Renming Song\thanks{Research supported in part by a grant from the Simons
Foundation (208236).} \hspace{1mm}\hspace{1mm} and \hspace{1mm}\hspace{1mm}
Rui Zhang
\hspace{1mm} }
\date{}
\maketitle

\bigskip\bigskip

{\narrower{\narrower

\centerline{\bf Abstract}

\bigskip
Suppose that $X=\{X_t, t\ge 0\}$ is a supercritical superprocess
on a locally compact separable metric space $(E, m)$.
Suppose that the spatial motion of $X$ is a Hunt process
satisfying certain conditions and that the branching mechanism is of
the form
$$
\psi(x,\lambda)=-a(x)\lambda+b(x)\lambda^2+\int_{(0,+\infty)}(e^{-\lambda y}-1+\lambda y)n(x,dy),
\quad x\in E, \quad\lambda> 0,
$$
where $a\in \mathcal{B}_b(E)$, $b\in \mathcal{B}_b^+(E)$ and $n$ is a kernel
from $E$ to $(0,\infty)$ satisfying
$$
  \sup_{x\in E}\int_0^\infty y^2 n(x,dy)<\infty.
$$
Put $T_tf(x)=\P_{\delta_x}\langle f,X_t\rangle$. Let $\lambda_0>0$ be the
largest eigenvalue of the generator $L$ of $T_t$,
and $\phi_0$ and $\wh{\phi}_0$ be the eigenfunctions of $L$ and $\widehat{L}$ (the dural of $L$)
respectively associated with $\lambda_0$.
Under some conditions on the spatial motion and the $\phi_0$-transformed semigroup of $T_t$,
we prove that for a large class of suitable functions $f$,
we have
$$
\lim_{t\rightarrow\infty}e^{-\lambda_0 t}\langle f, X_t\rangle =
W_\infty\int_E\wh{\phi}_0(y)f(y)m(dy),\quad \P_{\mu}\mbox{-a.s.},
$$
for any finite initial measure $\mu$ on $E$ with compact support,
where $W_\infty$ is the  martingale limit defined by
$W_\infty:=\lim_{t\to\infty}e^{-\lambda_0t}\langle \phi_0, X_t\rangle$.
Moreover, the exceptional set in the above limit does not depend on the initial
measure $\mu$ and the function $f$.
\bigskip


 \medskip


\par}\par}

\bigskip\bigskip

\begin{doublespace}

\section{Introduction}

Recently there have been quite a few papers on  law of large numbers for superdiffusions.
In \cite{E,ET,EW} some weak laws of large numbers (convergence
in law or in probability) were established.
The strong law of large numbers for superprocesses
was first studied
in \cite{CRW} followed by \cite{EKW14, KR, LRS, W}.
The continuity of the sample paths of the spatial motions played an important role
in all the papers mentioned above except \cite{CRW, KR}.
It is more difficult to establish strong law of large numbers for
superprocesses with discontinuous spatial motions.
For a good survey on recent developments in laws of large
numbers for branching Markov processes and superprocesses, see \cite{EKW14}.
In the papers mentioned above, either the spatial motion is assumed to be a diffusion,
or the spatial motion is assumed to be a symmetric Hunt process.
In the paper \cite{CRW} where the spatial motion is a symmetric Hunt process,
a condition on the smallness at ``infinity'' of the linear term in the branching mechanism
of the superprocess has to be assumed.
The purpose of this paper is to give a different setup under which the
strong law of large number for superprocesses holds.
The setup of this paper complements the previous setups. In particular,
the spatial motion may be discontinuous and non-symmetric.
We will give some examples satisfying the conditions of this paper.

The papers \cite{CRW, EKW14, LRS}
 dealt with strong law of large numbers for
superprocesses with spatially dependent branching mechansim.
The main ideas of the arguments of
\cite{CRW, EKW14, LRS}
are similar and consist of two steps.
The first step is to prove an almost sure limit result for discrete times,
and the second step is to prove that the result is  true for continuous times.
An essential difficulty comes from the second step.
\cite{LRS} gave a method for
the transition from lattice times to continuous times based on the resolvent
operator and approximation of the indicator function of an open subset
of $E$ by  resolvent functions.
The reason that this approximation works for superdiffusions
is that the sample paths of the spatial motion are continuous.
\cite{EKW14} also used this idea to show that indicator functions
can be approximated by resolvent functions.
For general superprocesses with spatial motions which might be discontinuous,
\cite{CRW} is
the first paper to establish a strong law of large numbers under a
second moment condition.
The paper \cite{CRW}
managed to overcome the difficulty of transition from discrete
times to continuous times with a highly non-trivial application of the
martingale formulation of superprocesses.
However, the assumptions of \cite{CRW} are restrictive in two aspects:
the spatial motion is assumed to be symmetric and the linear term
of the branching mechanism is assumed to satisfy a Kato class
condition at ``infinity''.

The papers \cite{W, KR} dealt with strong law of large numbers
for super-Brownian motions and super-$\alpha$-stable processes
with spatially independent branching mechanism respectively.
The key ingredients in the argument of \cite{W, KR}
are Fourier analysis and stochastic analysis, and the conditions
in \cite{W, KR} are quite different from those of
\cite{CRW, LRS}. The mean semigroup of the superprocess is assumed to
have a spectral gap in \cite{CRW, LRS}, while the mean semigroups of
the superprocesses of \cite{W, KR} have continuous spectra.
In this paper we assume that the spatial motion has a dual with respect
to a certain measure and that the branching mechanism satisfies
a second moment condition. Under the conditions of this paper, the mean
semigroup of the superprocess automatically has a spectral gap.

\subsection{Spatial process}\label{subs:sp}

Our assumptions on the underlying spatial process are the similar to those  in \cite{RSZ4}.
In this subsection, we recall the assumptions on the spatial process.

Suppose $(E,m)$ is a locally compact separable metric space and $m$ is
a $\sigma$-finite Borel measure
on $E$ with full support.
Let $E_\partial = E\cup\{\partial\}$ be the one-point compactification of $E$.
Every function $f$ on $E$ is automatically extended to
$E_\partial$
by setting $f(\partial)=0$.
 We will assume that $\xi=\{\xi_t,\Pi_x\}$ is a Hunt process on $E$ and $\zeta:=
 \inf\{t>0: \xi_t=\partial\}$ is the lifetime of $\xi$.
The transition semigroup  of $\xi$ will be denoted by $\{P_t,t\geq 0\}$.
We will always assume
that there exists a family of strictly positive
continuous functions
$\{p(t,x,y), t>0\}$ on $E\times E$ such that
$$  P_tf(x)=\int_E p(t,x,y)f(y)\,m(dy).$$
Define
\begin{equation}\label{e:RS1}
a_t(x):=\int_E p(t,x,y)^2\,m(dy),\qquad \widehat{a}_t(x):=\int_E p(t,y,x)^2\,m(dy).
\end{equation}
In this paper, we assume that

\begin{assumption}\label{assum1}
\begin{description}
  \item[(a)]
  For all $t>0$ and  $x\in E$, $\int_E p(t,y,x)\,m(dy)\le 1$.

  \item[(b)] For any $t>0$,
  $a_t$ and $\widehat{a}_t$ are continuous
     $L^1(E;m)$-integrable functions.

  \item[(c)] There exists $t_0>0$ such that $a_{t_0},\widehat{a}_{t_0}\in L^2(E;m)$.
\end{description}
\end{assumption}
By the Chapman-Kolmogorov equation and the Cauchy-Schwarz inequality,
\
\begin{equation}\label{1.11}
  p(t+s,x,y)=\int_E p(t,x,z)p(s,z,y)\,m(dz)\le (a_t(x))^{1/2}(\widehat{a}_s(y))^{1/2}.
\end{equation}
Therefore,
$a_{t+s}(x)\le \int_E \widehat{a}_s(y)\,m(dy)a_t(x)$ and $\widehat{a}_{t+s}(x)\le
\int_E a_s(y)\,m(dy)\widehat{a}_t(x).$
Thus under condition (b), the condition (c)
above is equivalent to
 \begin{description}
   \item[(c$'$)] There exists $t_0>0$ such that for all $t\ge t_0$, $a_{t},\widehat{a}_t\in L^2(E;m)$.
 \end{description}

Under Assumption \ref{assum1}(a), for every $t>0$, both $P_t$ and
the operator $\widehat P_t$ defined by
$\widehat P_t f(x)= \int_E p(t, y, x) f(y) m(dy)$
are contraction operators in $L^p(E; m)$ for every $p\in [1, \infty]$, and they are dual
to each other.
Assumption \ref{assum1}(b) implies
that each
$P_t$ is a Hilbert-Schmidt operator in $L^2(E; m)$
and thus is compact. Hence $P_t$  has discrete spectrum.

\subsection{Superprocesses}\label{superp}

In this subsection, we introduce the superprocesses.
Let $\mathcal{B}_b(E)$ (respectively, $\mathcal{B}_b^+(E)$) be the family of
bounded (respectively, nonnegative bounded) Borel functions on $E$.
Denote by $\langle\cdot,\cdot\rangle_m$ the inner product in $L^2(E;m)$.

The superprocess $X=\{X_t, t\ge 0\}$
is determined by three parameters:
a spatial motion $\xi=\{\xi_t, \Pi_x\}$ on $E$ satisfying the
assumptions of the previous subsection,
a branching rate function $\beta(x)$ on $E$ which is a nonnegative bounded Borel
function and a branching mechanism $\psi$ of the form
\begin{equation}
\psi(x,\lambda)=-a(x)\lambda+b(x)\lambda^2+\int_{(0,+\infty)}(e^{-\lambda y}-1+\lambda y)n(x,dy),
\quad x\in E, \quad\lambda> 0,
\end{equation}
where $a\in \mathcal{B}_b(E)$, $b\in \mathcal{B}_b^+(E)$ and $n$ is a kernel
from $E$ to $(0,\infty)$ satisfying
\begin{equation}\label{n:condition}
  \sup_{x\in E}\int_0^\infty y^2 n(x,dy)<\infty.
\end{equation}

Let ${\cal M}_F(E)$ be the space of finite measures on $E$,
equipped with the weak convergence  topology.
As usual, $\langle f,\mu\rangle:=\int f(x)\mu(dx)$ and $\|\mu\|:=\langle 1,\mu\rangle$.
According to \cite[Theorem 5.12]{Li11}, there is a Borel right process
$X=\{\Omega, {\cal G}, {\cal G}_t, X_t, \P_\mu\}$ taking values in
$\mathcal{M}_F(E)$, called superprocess,  such that
for every
$f\in \mathcal{B}^+_b(E)$ and $\mu \in \mathcal{M}_F(E)$,
\begin{equation}
  -\log \P_\mu\left(e^{-\langle f,X_t\rangle}\right)=\langle u_f(\cdot,t),\mu\rangle,
\end{equation}
where $u_f(x,t)$ is the unique positive solution to the equation
\begin{equation}\label{1.3}
  u_f(x,t)+\Pi_x\int_0^t\psi(\xi_s, u_f(\xi_s,t-s))\beta(\xi_s)ds=\Pi_x f(\xi_t),
\end{equation}
where $\psi(\partial,\lambda)=0, \lambda>0$.
Here
$({\cal G}, {\cal G}_t)_{t\ge 0}$ are augmented,
$(\mathcal{G}_t, t\ge 0)$ is right continuous
and $X$ satisfies the Markov property with respect to $(\mathcal{G}_t, t\ge 0)$.
Moreover, such a superprocess $X$ has a Hunt realization in ${\cal M}_F(E)$,
see \cite[Theorem 5.12]{Li11}.
In this paper, the superprocess we deal with
always takes such a  Hunt realization.

Define
\begin{equation}\label{e:alpha}
\alpha(x):=\beta(x)a(x)\quad \mbox{and }
A(x):=\beta(x)\left( 2b(x)+\int_0^\infty y^2 n(x,dy)\right).
\end{equation}
Then, by our assumptions, $\alpha(x)\in \mathcal{B}_b(E)$ and $A(x)\in\mathcal{B}_b^+(E)$.
Thus there exists $K>0$ such that
\begin{equation}\label{1.5}
  \sup_{x\in E}\left(|\alpha(x)|+A(x)\right)\le K.
\end{equation}
For any $f\in\mathcal{B}_b(E)$ and $(t, x)\in (0, \infty)\times E$, define
\begin{equation}\label{1.26}
   T_tf(x):=\Pi_x \left[e^{\int_0^t\alpha(\xi_s)\,ds}f(\xi_t)\right].
\end{equation}
It is well-known that $T_tf(x)=\P_{\delta_x}\langle f,X_t\rangle$ for every $x\in E$.
It is known that (see, e.g., \cite{RSZ4} and \cite[Lemma 2.1]{RSZ6})
$\{T_t,t\ge 0\}$ is a strongly
continuous semigroup on $L^2(E;m)$ and there exists a function
$q(t, x, y)$ on $(0, \infty)\times E\times E$ which is continuous in $(x, y)$ for each $t>0$
such that
\begin{equation}\label{comp}
 e^{-Kt}p(t,x,y) \le q(t,x,y)\le e^{Kt}p(t,x,y)  \quad \hbox{for }
 (t, x, y)\in (0, \infty)\times E\times E
\end{equation}
and that for any bounded Borel function $f$ on $E$ and $(t, x)\in (0, \infty)\times E$,
$$
T_tf(x)=\int_Eq(t, x, y)f(y)m(dy).
$$
Define
\begin{equation}\label{e:RS2}
b_t(x):=\int_E q(t,x,y)^2\,m(dy),\qquad \widehat{b}_t(x):=\int_E q(t,y,x)^2\,m(dy).
\end{equation}
Then $b_t$ and
$\widehat{b}_t$ enjoy the following properties:
\begin{description}
  \item[(i)] For any $t>0$, we have $b_t,  \widehat{b}_t\in L^1(E;m)$. Moreover, $b_t(x)$ and
  $\widehat{b}_t(x)$ are continuous
  in $x\in E$.

   \item[(ii)] There exists $t_0>0$ such that for all $t\ge t_0$, $b_{t},\widehat{b}_{t}\in L^2(E;m)$.
\end{description}
 Let $\{\widehat{T}_t, t>0\}$ be the adjoint semigroup of $\{T_t, t\ge 0\}$
on $L^2(E, m)$ defined by
$$
\widehat{T}_tg(x)=\int_E q(t,y,x)g(y)\,m(dy).
$$
It is easy to see $\widehat T_t$ is the dual operator of $T_t$ in $L^2(E; m)$. It follows
that $\{\widehat{T}_t, t>0\}$ is also strongly continuous
in $L^2(E, m)$. Since $q(t,\cdot,y)$  and $a_{t}$ are continuous, by \eqref{1.11} and \eqref{comp},
using the dominated convergence theorem, we get that  for any $t>0$ and $f\in L^2(E; m)$, $T_tf$ and $\widehat{T}_tf$
are continuous.

It follows from (i) above that,
for any $t>0$, $T_t$ and $\widehat{T}_t$ are compact operators in $L^2(E; m)$.
Let $L$ and $\widehat{L}$ be the infinitesimal generators of the semigroups
$\{T_t\}$ and $\{\widehat{T}_t\}$ in $L^2(E;m)$ respectively.
Let $\sigma(L)$ and $\sigma(\widehat{L})$ be the spectra of
$L$ and $\widehat{L}$. It follows from \cite[Theorem 2.2.4 and Corollary 2.3.7]{Pazy}
that both $\sigma(L)$ and $\sigma(\widehat{L})$ consist of eigenvalues, and that
$\sigma(L)$ and $\sigma(\widehat{L})$ have the same number, say $N$, of eigenvalues.
Let $\mathbb{I}=\{0, \cdots, N-1\}$ if $N<\infty$ and $\mathbb{I}=\{0, \cdots\}$ otherwise.
Define $\lambda_0:=\sup \Re(\sigma(L))=\sup\Re(\sigma(\widehat{L}))$.
By Jentzsch's theorem (Theorem V.6.6 on page 337 of \cite{Sch}),
$\lambda_0$ is an eigenvalue of multiplicity 1 for both $L$ and $\widehat{L}$.
Assume that $\phi_0$ and
${\widehat \phi}_0$
are the eigenfunctions of $L$ and $\widehat{L}$
respectively associated with $\lambda_0$.
$\phi_0$ and ${\widehat \phi}_0$ can be chosen to be continuous strictly
positive and satisfy $\|\phi_0\|_2=1$ and $\langle \phi_0,{\widehat \phi}_0\rangle_m=1$.
We list the eigenvalues of $\{\lambda_k, k\in\mathbb{I}\}$ of $L$
in an order so that $\lambda_0>\Re(\lambda_1)\ge \Re(\lambda_2)\ge\cdots$.
Then $\{\overline{\lambda}_k, k\in\mathbb{I}\}$ are the eigenvalues of $\widehat{L}$.
For convenience, we define, for any positive integer not in $\mathbb{I}$,
$\lambda_k=\overline{\lambda}_k=-\infty$. For $k\in \mathbb{I}$, we write
$\Re_k:=\Re(\lambda_k)$.
We use the convention $\Re_\infty=-\infty$.

For $t>0$, $T_t\phi_0(x)=e^{\lambda_0t}\phi_0(x)$, and thus
\begin{equation}\label{1.1}
  \phi_0(x)\le e^{-\lambda_0t}b_{t}(x)^{1/2}.
\end{equation}
Similarly, we have $\widehat{T}_t{\widehat \phi}_0(x)=e^{\lambda_0t}{\widehat \phi}_0(x)$ and ${\widehat \phi}_0(x)\le e^{-\lambda_0t}\|{\widehat \phi}_0\|_2\widehat{b}_{t}(x)^{1/2}.$
Therefore, by Assumption \ref{assum1}(c),
$\phi_0\in L^2(E;m)\cap L^4(E;m)$.
In this paper, we always assume that the superprocess $X$ is supercritical, that is, $\lambda_0>0$.
Define $W_t:=e^{-\lambda_0t}\langle \phi_0,X_t\rangle$. By the Markov
property of $X$, $\{W_t,t\ge 0\}$ is a nonnegative
martingale with respect to $\{\mathcal{G}_t,t\ge0\}$, and thus the
 $W_\infty:=\lim_{t\to\infty}W_t$ exists. Under our assumptions,
 $W_t$ is a $L^2$-bounded martingale, thus $W_\infty$ is non-degenerate, that is $\P_{\mu}(W_\infty>0)>0$.

\subsection{Main results}

In this subsection, we state our main results. In the remainder of this paper,
whenever we talk about an initial configuration $\mu\in {\cal M}_F(E)$, we
always implicitly assume that it has compact support.

For $q>\max\{K, \lambda_0\}$ and $f\in L^p(E;m)$ with $p\geq 1$, define,
$$
U_qf(x):=\left\{
             \begin{array}{ll}
               \displaystyle\int_0^\infty e^{-qs}T_sf(x)\,ds, &\mbox{if }\displaystyle\int_0^\infty
               e^{-qs}T_s|f|(x)\,ds<\infty; \\
              0,
              &\displaystyle \hbox{otherwise.}
             \end{array}
           \right.
$$
Note that for $p \geq 1$, by Assumption \ref{assum1}(a)  and \eqref{comp}
\begin{eqnarray}\label{L_p}
 && \left(\int_E \Big(\int_0^\infty e^{-qs}T_s|f|(x)\,ds\Big)^p\,m(dx)\right)^{1/p}\nonumber\\
 & \le& \int_0^\infty e^{-qs}\|T_s(|f|)\|_p\,ds\nonumber\\
  &\le& \int_0^\infty e^{-qs}e^{Ks}\,ds\|f\|_p<\infty,
\end{eqnarray}
which implies that $\int_0^\infty e^{-qs}T_s|f|(x)\,ds \in L^p(E;m)$,
and thus $\int_0^\infty e^{-qs}T_s|f|(x)\,ds<\infty, m$-a.e.
Consequently, $U_qf\in L^p(E;m)$.
In Lemma \ref{lem2} below,
we will show that if $f\in L^2(E,m)\cap L^4(E,m)$ then $\langle U_qf,X_t\rangle$ is well defined.

\begin{thrm}
\label{thrm1}
Assume that Assumption \ref{assum1} holds.
If $g=U_qf $ for some $f\in L^2(E, m)\cap L^4(E;m)$ and $q>\max\{K, \lambda_0\}$,
then for any $\mu\in \mathcal{M}_F(E)$, as $t\to\infty$,
\begin{equation}\label{7.35}
  e^{-\lambda_0t}\langle g, X_t\rangle\to \langle g,{\widehat \phi}_0\rangle_mW_\infty,
   \quad  \P_{\mu}\mbox{-a.s.}
\end{equation}
\end{thrm}
For any $f\ge 0$, define
\begin{equation}\label{h-transf}T^{\phi_0}_tf(x)=
\frac{e^{-\lambda_0t}}{\phi_0(x)}\Pi_x\left[\exp\left(\int^t_0\alpha(\xi_s)ds\right)
(f\phi_0)(\xi_t)\right].\end{equation}

Let $C_0(E;\R)$ denote the family of real-valued continuous functions $f$ on $E$ with the property
that $\lim_{x\to\partial}f(x)=0$.

We will also make the following assumption in this paper.

\begin{assumption}\label{assum2} The semigroup $\{T^{\phi_0}_t, t\ge 0\}$ has the following properties:
For any $f\in C_0(E;\R)$,
\begin{equation}\label{Feller}\lim_{t\to 0}\|T^{\phi_0}_tf-f\|_\infty=0.\end{equation}
\end{assumption}

The following theorem is the main result of this paper.

\begin{thrm}\label{main theorem}
Under Assumptions \ref{assum1} and \ref{assum2}, there
exists $\Omega_0\subset\Omega$ of probability one (that is,
$\mathbb P_\mu (\Omega_0)=1$ for every
$\mu\in \mathcal M_F(E)$) such that,
for every $\omega\in\Omega_0$ and for every bounded Borel
function $f$ on $E$ satisfying (a) $|f|\le c\phi_0$ for some $c>0$ and (b)
the set of  discontinuous points of $f$ has zero $m$-measure, we have
\begin{eqnarray}\label{main result2}
\lim_{t\rightarrow\infty}e^{-\lambda_0 t}\langle f, X_t\rangle (\omega) =
W_\infty(\omega)\int_E{\widehat \phi}_0(y)f(y)m(dy).
\end{eqnarray}
\end{thrm}

Assumption \ref{assum2} will be used to extend the test functions from
resolvent functions $g=U_qf $ with $f\in L^2(E, m)\cap L^4(E;m)$
to functions of the form $g=f\phi_0$ with $f\in C_0(E;\R)$.
We will give some examples in Section \ref{s:example} to show that
Assumptions \ref{assum1} and \ref{assum2} are satisfied
by many interesting superprocesses including
super Ornstein-Uhlenbeck processes
(both inward and outward)
 and superprocesses with discontinuous
spatial motions.

\begin{remark}{\rm

(1) Compared with \cite{CRW}, our spatial motion
can be  nonsymmetric and we do not assume that $\alpha(x)=\beta (x) a(x)$ is in the Kato class $K_\infty(\xi)$. The latter would require $\alpha$ be
 in some sense small at $\infty$
(see \cite{CRW} for the definition of  $K_\infty(\xi)$). In \cite{CRW},
a compact  embedding condition (see \cite[2.4]{CRW}) is also assumed
to ensure that the generator of the semigroup $\{T_t,t\ge 0\}$  has a spectral gap.
In this paper, we assume instead
Assumption \ref{assum1}, which implies
that the generator of $\{T_t, t\ge 0\}$ has discrete spectrum.

(2) Compared with \cite{LRS} where the spatial motion is a diffusion,
our spatial motion may be discontinuous.
The setup of \cite{LRS} and the setup of the present are also different in the
following ways.
In \cite{LRS}, the semigroup of the spatial motion is assumed
to be intrinsic ultracontractive.
This condition is pretty strong and it excludes
some interesting examples including the OU process.
In this paper, we assume Assumption \ref{assum1} instead, which is weaker than
the intrinsic ultracontractive property and is enough to insure that, for resolvent functions $g$,
the limit $\lim_{t\to\infty}e^{-\lambda_0 t}\langle g, X_t\rangle$ exists almost surely.
In \cite{LRS}, the branching mechanism is assumed to satisfy a $L\log L$ condition,
while in this paper, we assume that the branching mechanism satisfies
a second moment condition.
}
\end{remark}

\section{Preliminaries}

\subsection{Moment estimates}

By \cite[Lemma 2.2]{RSZ4} with $k=1$, for any $t_1>0$
and $a<-\Re(\lambda_1)$, there exists a constant $c=c(a, t_1)>0$ such that
for all $(t, x, y)\in (2t_1, \infty)\times E\times E$,
\begin{equation}\label{e:RSZ42.3}
\left|q(t, x, y)-e^{\lambda_0t}\phi_0(x){\widehat \phi}_0(y) \right|
\le ce^{-at}b_{t_1}(x)^{1/2}\widehat{b}_{t_1}(y)^{1/2}.
\end{equation}
Multiplying both sides by $e^{-\lambda_0t}$, we get that
for all $(t, x, y)\in (2t_1, \infty)\times E\times E$,
$$
\left|e^{-\lambda_0t}q(t,x,y)-\phi_0(x){\widehat \phi}_0(y)\right|
   \le ce^{-(a+\lambda_0)t}b_{t_1}(x)^{1/2}\widehat{b}_{t_1}(y)^{1/2}.
$$
Note that $a<-\Re(\lambda_1)$ is equivalent to $a+\lambda_0<\lambda_0-\Re(\lambda_1)$.
Thus for any $\widetilde{a}\in (0, \lambda_0-\Re(\lambda_1))$ and
$t_1>0$,  there exists $c_1=c_1(\widetilde{a}, t_1)>0$ such that
for all $(t, x, y)\in (2t_1, \infty)\times E\times E$,
\begin{equation}\label{density}
  \left|e^{-\lambda_0t}q(t,x,y)-\phi_0(x){\widehat \phi}_0(y)\right|
   \le c_1e^{-\widetilde{a}t}b_{t_1}(x)^{1/2}\widehat{b}_{t_1}(y)^{1/2}.
\end{equation}
Thus, for $f\in L^2(E;m)$, we have for all
$(t, x)\in (2t_1, \infty)\times E$,
$$
\left|e^{-\lambda_0t}T_tf(x)-\phi_0(x)\langle f,{\widehat \phi}_0\rangle_m\right|
   \le c_1\|\widehat{b}_{t_1}^{1/2}\|_2\|f\|_2e^{-\widetilde{a}t}b_{t_1}(x)^{1/2},
$$
which implies that there exists $c_2=c_2(\widetilde{a}, t_1)>0$ such that
for all $(t, x)\in (2t_1, \infty)\times E$,
\begin{equation}\label{moment1}
  \left|e^{-\lambda_0t}T_tf(x)-\phi_0(x)\langle f,{\widehat \phi}_0\rangle_m\right|
   \le c_2\|f\|_2e^{-\widetilde{a}t}b_{t_1}(x)^{1/2}.
\end{equation}
Hence, by \eqref{1.1}, we have
\begin{eqnarray*}
e^{-\lambda_0t}|T_tf(x)|&\le& \phi_0(x)|\langle f,{\widehat \phi}_0\rangle_m|+c_2\|f\|_2
e^{-\tilde{a}t}b_{t_1}(x)^{1/2}\\
&\le& (e^{-\lambda_0t_1}\|{\widehat \phi}_0\|_2+c_2)\|f\|_2b_{t_1}(x)^{1/2}.
\end{eqnarray*}
Thus there exists $c_3=c_3(\widetilde{a}, t_1)>0$ such that
for all $(t, x)\in (2t_1, \infty)\times E$,
\begin{equation}\label{moment3}
  |T_tf(x)|\le c_3\|f\|_2e^{\lambda_0t}b_{t_1}(x)^{1/2}.
\end{equation}

We now recall the second moment formula for the superprocess $\{X_t, t\ge 0\}$
(see, for example, \cite{RSZ3}):
for $f\in L^2(E;m)\cap L^4(E;m)$ and $\mu\in\mathcal{M}_F(E)$, we have for any $t >0$,
\begin{equation}\label{1.13}
   {\V}{\rm ar}_{\mu}\langle f,X_t\rangle=\langle{\V}{\rm ar}_{\delta_\cdot}
   \langle f,X_t\rangle, \mu\rangle
   =\int_E\int_0^tT_{s}[A(T_{t-s}f)^2](x)\,ds\mu(dx),
\end{equation}
where $\mathbb{V}{\rm ar}_{\mu}$ stands for the variance under $\P_{\mu}$
and $A(x)$ is the function defined in \eqref{e:alpha}.
Moreover, for $f\in L^2(E;m)\cap L^4(E;m)$,
\begin{equation}\label{1.4}
  \V{\rm ar}_{\delta_x}\langle f,X_t\rangle\le e^{Kt}T_t(f^2)(x)\in L^2(E;m).
\end{equation}

In the following lemma, we give a useful estimate on the second moment of $X$.
If we choose the constant $\widetilde{a}\in (0, \lambda_0- \Re(\lambda_1 ))$ small enough,
 we can get the next lemma by \cite[Lemma 2.5]{RSZ4}.
Here we give a direct proof.

\begin{lemma}\label{lem1}
Suppose that Assumption \ref{assum1} holds.
 For any $\widetilde{a}\in (0, (\lambda_0-\Re(\lambda_1))\wedge(\lambda_0/2))$
and $f\in L^2(E;m)\cap L^4(E;m)$ with $\langle f,{\widehat \phi}_0\rangle_m=0$,
there exists $c_4=c_4(t_0, \widetilde{a}, f)>0$ such that
\begin{equation}\label{moment2}
  \sup_{t>10t_0}e^{2(-\lambda_0+\widetilde{a})t}\V{\rm ar}_{\delta_x}\langle f,
  X_t\rangle\le c_4b_{t_0}(x)^{1/2}.
\end{equation}
\end{lemma}

\textbf{Proof:}
In the following proof, we use $c=c(t_0, \tilde a, f)$ to denote a constant whose value may change from one appearance to another.
Recall that
$$\V{\rm ar}_{\delta_x}\langle f,X_t\rangle=\Big(\int_0^{2t_0}+
\int_{2t_0}^{t-2t_0}+\int_{t-2t_0}^{t}\Big)
T_{s}[A(T_{t-s}f)^2](x)\,ds.$$
In the following we will deal with the above three parts separately.

(i) For $t>10t_0$ and $s<2t_0$, by \eqref{moment1}, we have
$$
|T_{t-s}f(x)|
\le c
e^{(\lambda_0-\widetilde{a})(t-s)}b_{4t_0}(x)^{1/2}.$$
Thus,
\begin{equation*}
  \int_0^{2t_0}T_{s}[A(T_{t-s}f)^2](x)\,ds
  \le c
  e^{2(\lambda_0-\widetilde{a})t}
  \int_0^{2t_0}T_{s}[b_{4t_0}](x)\,ds.
\end{equation*}
If we can prove that
\begin{equation}\label{toshow}\int_0^{2t_0}T_{s}[b_{4t_0}](x)\,ds\le c b_{t_0}(x)^{1/2},\end{equation}
we will get
\begin{equation}\label{2.5}
  \int_0^{2t_0}T_{s}[A(T_{t-s}f)^2](x)\,ds
  \le c
  e^{2(\lambda_0-\widetilde{a})t}b_{t_0}(x)^{1/2}.
\end{equation}
Now we prove \eqref{toshow}. By Fubini's theorem and  H\"{o}lder's inequality, we get
\begin{eqnarray*}
  a_{t+s}(x) &=& \int_E p(t+s,x,y)\int_{E}p(t,x,z)p(s,z,y)\,m(dz)\,m(dy)\\
&=& \int_E p(t,x,z)\int_{E}p(t+s,x,y)p(s,z,y)\,m(dy)\,m(dz)\\
&\le&  a_{t+s}(x)^{1/2}\int_E p(t,x,z)a_s(z)^{1/2}\,m(dz)
\end{eqnarray*}
which implies
\begin{equation}\label{8.9}
  a_{t+s}(x)\le \left(\int_E p(t,x,z)a_s(z)^{1/2}\,m(dz)\right)^2\le \int_E p(t,x,z)a_s(z)\,m(dz).
\end{equation}By \eqref{8.9}, we get
$$
b_{4t_0}(x)\le e^{8Kt_0}a_{4t_0}(x)\le e^{10Kt_0}T_{2t_0}(a_{2t_0})(x).
$$
Thus, by Assumption \ref{assum1}(c$'$) and \eqref{moment3}, we have
\begin{eqnarray}\label{1.37}
  \int_{0}^{2t_0}T_s(b_{4t_0})(x)\,ds&\le& e^{10Kt_0}\int_{0}^{2t_0}T_{s+2t_0}(a_{2t_0})(x)\,ds\nonumber\\
  &\le& c\int_{0}^{2t_0}e^{\lambda_0(s+2t_0)}\,ds b_{t_0}(x)^{1/2}\le c b_{t_0}(x)^{1/2}.
\end{eqnarray}
Therefore \eqref{toshow} holds.

(ii) For $t>10t_0$ and $s\in(2t_0,t-2t_0)$,
by \eqref{moment1}, \eqref{moment3} and  Assumption \ref{assum1}(c$'$),
$$
T_{s}[A(T_{t-s}f)^2](x)
 \le c
 e^{2(\lambda_0-\widetilde{a})(t-s)}T_s(b_{t_0})(x)
 \le c
  e^{2(\lambda_0-\widetilde{a})(t-s)}e^{\lambda_0s}b_{t_0}(x)^{1/2}.
$$
Thus, using the fact $\lambda_0-2\widetilde{a}>0$,
\begin{eqnarray}\label{3.2}
  \int_{2t_0}^{t-2t_0}T_{s}[A(T_{t-s}f)^2](x)\,ds
  &\le& c
  e^{2(\lambda_0-\widetilde{a})t}\int_{2t_0}^{t-2t_0}e^{-(\lambda_0-2\widetilde{a})s}\,dsb_{t_0}(x)^{1/2}\nonumber\\
  &\le& c
  e^{2(\lambda_0-\widetilde{a})t}b_{t_0}(x)^{1/2}.
\end{eqnarray}

(iii) For $t>10t_0$ and $s>t-2t_0$, since $|T_{t-s}f(x)|^2\le e^{K(t-s)}T_{t-s}(f^2)(x)$,
\begin{eqnarray*}
T_{s}[A(T_{t-s}f)^2](x)&\le& Ke^{K(t-s)}T_t(f^2)(x)\le Ke^{2t_0K} c_3e^{\lambda_0t}b_{t_0}(x)^{1/2}\\
&\le& Ke^{2t_0K} c_3e^{2(\lambda_0-\widetilde{a})t}b_{t_0}(x)^{1/2},
\end{eqnarray*}
where in the last equality we use the fact $\lambda_0-2\widetilde{a}>0$.
Thus,
\begin{equation}\label{3.3}
  \int_{t-2t_0}^tT_{s}[A(T_{t-s}f)^2](x)\,ds
  \le c
  e^{2(\lambda_0-\widetilde{a})t}b_{t_0}(x)^{1/2}.
\end{equation}
Combining \eqref{2.5}, \eqref{3.2} and \eqref{3.3}, we get \eqref{moment2}.\hfill$\Box$

\subsection{Martingale measure for superprocesses}

In this subsection, we recall the
associated martingale measure
for the superprocess $X$.
For more details, see, for instance, \cite[Chapter 7]{Li11}.
The martingale measure
for superprocesses is a very useful tool in the proof of our main theorems.

For our superprocess $X$,
there exists a worthy $(\mathcal{G}_t)$-martingale measure
$\{M_t(B)=M(t, B); t\ge 0, B\in\mathcal{B}(E)\}$
with covariation measure
$$
  \nu(ds,dx,dy):=ds\int_EA(z)\delta_z(dx)\delta_z(dy)X_s(dz)
$$
such that for $t\ge 0$ and $f\in L^2(E;m)\cap L^4(E;m)$, we have, $\P_{\mu}$-a.s.,
\begin{equation}\label{5.1}
  \langle f,X_t\rangle=\langle T_tf,\mu\rangle+\int_0^t\int_E T_{t-s}f(z)\,M(ds,dz).
\end{equation}
For any $u>0$ and
$0\le t\le u$,
we define
$$M^{(u)}_t:=\int_{0}^{t}\int_E T_{u-s}f(x)M(ds,dx).$$
Then, for any $\mu\in \mathcal{M}_F(E)$,
$\{M^{(u)}_t,0\le t\le u\}$ is a
cadlag square-integrable martingale under $\P_\mu$ with
\begin{equation}\label{5.22}
 \langle M^{u}\rangle_t=\int_0^t \langle A(T_{u-s}f)^2, X_s\rangle\,ds.
\end{equation}
Here cadlag means ``right continuous having left limits".
Note that
\begin{equation}\label{6.10}
  \P_{\mu}(M^{(u)}_u)^2=\P_{\mu}\langle M^{u}\rangle_u=\V ar_{\mu}\langle f,X_u\rangle.
\end{equation}

In the remainder of this paper, we will always assume that
$q>\max\{K, \lambda_0\}$.

\begin{lemma}\label{lem2}
Assume that Assumption \ref{assum1} holds. If $f\in L^2(E;m)\cap L^4(E;m)$,
then for any $\mu\in {\cal M}_F(E)$,
$$
\P_{\mu}\left(\langle U_q|f|,X_t\rangle<\infty \hbox{ for } t\ge0\right
)=\P_{\mu}\left(\langle U_qf,X_t\rangle \mbox{ is finite}
\hbox{ for } t\ge0\right)=1.
$$
Moreover, $\P_{\mu}$-a.s., $\langle U_qf,X_t\rangle$ is
cadlag on $[0,\infty)$,
and for all $t>0$,
\begin{eqnarray}\label{5.6}
  &&\langle U_qf,X_{t}\rangle
  =\langle T_{t}(U_qf),\mu\rangle
  +e^{qt} \int_{t}^\infty e^{-qu}M^{(u)}_t\,du.
\end{eqnarray}
\end{lemma}

\noindent\textbf{Proof:}
When the spatial motion $\xi$ is symmetric, this lemma has been established in
\cite[lemma 2.4 and Lemma 2.5]{RSZ5}. The proof for the non-symmetric case is almost the same.
For reader's convenience, we include a proof here.
We can check that the argument in the proof of \cite[Lemma 2.4]{RSZ5} works
without the assumption that $\xi$ is $m$-symmetric,
so $\langle U_qf,X_t\rangle$ is
right continuous on $[0,\infty)$, $\P_{\mu}$-a.s.

For $f\in L^2(E;m)\cap L^4(E;m)$, $U_qf\in L^2(E;m)\cap L^4(E;m)$.
By \eqref{5.1}, for $t>0$ and $\mu\in \mathcal{M}_F(E)$, we have, $\P_{\mu}$-a.s.,
\begin{eqnarray}
  \langle U_qf,X_{t}\rangle
  &=&\langle T_{t}(U_qf),\mu\rangle+\int_0^{t}\int_E T_{t-s}(U_qf)(z)M(ds,dz)\nonumber\\
  &=& \langle T_{t}(U_qf),\mu\rangle
  +\int_0^{t}\int_E \int_{0}^\infty e^{-qu}T_{u+t-s}f(z)\,duM(ds,dz)\nonumber\\
  &=& \langle T_{t}(U_qf),\mu\rangle
  +e^{qt}\int_0^{t}\int_E \int_{t}^\infty e^{-qu}T_{u-s}f(z)\,duM(ds,dz)\nonumber\\
  &=&\langle T_{t}(U_qf),\mu\rangle
  +e^{qt} \int_{t}^\infty e^{-qu}\,du\int_0^{t}\int_ET_{u-s}f(z)M(ds,dz)\nonumber\\
  &:=&J_1^f(t)+e^{qt}J_2^f(t),\label{5.8}
\end{eqnarray}
where the fourth equality follows from the stochastic Fubini's theorem for martingale measures
(see, for instance, \cite[Theorem 7.24]{Li11}).
Thus,  for $t>0$ and $\mu\in \mathcal{M}_F(E)$,
\begin{equation}\label{4.3}
  \P_{\mu}\left(\langle U_qf,X_t\rangle=J_1^f(t)+e^{qt}J^f_2(t)\right)=1.
\end{equation}
Then, in light of \eqref{4.3}, to prove \eqref{5.6},
it suffices to prove that $J_1^f(t)$ and $J_{2}^f(t)$ are all cadlag in $(0,\infty)$, $\P_{\mu}$-a.s.
For $J_1^f(t)$, by Fubini's theorem, for $t>0$,
$$J_1^f(t)=e^{qt}\int_{t}^\infty e^{-qs}\langle T_sf,\mu\rangle\,ds.$$
Thus, it is easy to see that $J_1^f(t)$ is continuous in $t\in (0,\infty)$.
Now, we consider $J_{2}^f(t)$.
We claim that, for any $t_1>0$,
\begin{equation}\label{4.27'}
  \P_{\mu}\left( J_{2}^f(t) \mbox{ is cadlag in } [t_1,\infty)\right)=1.
\end{equation}
By the definition of $J_{2}^f$, for $t\ge t_1 $,
$$
  J_{2}^f(t)=\int_{t_1}^{\infty}e^{-qu}M^{(u)}_{t}1_{t<u}\,du.
$$
Since $t\mapsto M^{(u)}_{t}\textbf{1}_{t<u}$ is right continuous,
by the dominated convergence theorem, to prove \eqref{4.27'}, it suffices to show that
\begin{equation}\label{4.8}
  \P_{\mu}\left(\int_{t_1}^{\infty}e^{-qu}\sup_{t\ge t_1}\left(|M^{(u)}_{t}|
  1_{t<u}\right)\,du<\infty\right)=1.
\end{equation}
By  the $L^p$-maximum inequality and \eqref{6.10}, we have
\begin{eqnarray}\label{4.10}
  &&\P_{\mu}\left(\int_{t_1}^{\infty}e^{-qu}\sup_{t\ge t_1}\left(|M^{(u)}_{t}|1_{t<u}\right)\,du\right)
  \le 2\int_{t_1}^\infty e^{-qu}\sqrt{\P_{\mu}\left|M^{(u)}_{u}\right|^2}\,du\nonumber\\
  &=&2\int_{t_1}^\infty e^{-qu}\sqrt{\int_E\V ar_{\delta_x}\langle f,X_u\rangle\,\mu(dx)}\,du.
\end{eqnarray}
By \eqref{1.4} and \eqref{moment3}, we have, for $u>t_1$,
\begin{eqnarray*}\int_E\V ar_{\delta_x}\langle f,X_u\rangle\,\mu(dx)&\le& e^{Ku}\int_ET_u(f^2)(x)\,\mu(dx)\\
&\le& c
e^{Ku}e^{\lambda_0u}\int_Eb_{t_1/2}(x)^{1/2}\,\mu(dx),\end{eqnarray*}
where $c=c(t_1, \tilde a, f)$ is a positive constant
and $b_t (x)$ is the function defined in \eqref{e:RS2}.
Since $x\mapsto b_{t_1/2}(x)$ is continuous and $\mu$ has compact support,
we have $\int_E b_{t_1/2}(x)^{1/2}\mu(dx)<\infty.$
Thus by \eqref{4.10}, we have
\begin{eqnarray*}
  &&\P_{\mu}\left(\int_{t_1}^{\infty}e^{-qu}\sup_{t\ge t_1}\left(|M^{(u)}_{t}|1_{t<u}\right)\,du\right)\\
  &\le& 2\sqrt{c}
   \int_{t_1}^\infty e^{-qu}e^{(K+\lambda_0)u/2}\,du\sqrt{\int_E b_{t_1/2}(x)^{1/2}\,\mu(dx)}<\infty.
\end{eqnarray*}
Now \eqref{4.8} follows immediately.
Since $t_1>0$ is arbitrary, we have
$$
  \P_{\mu}\left( J_{2}^f(t) \mbox{ is cadlag in } (0,\infty)\right)=1.
$$
The proof is now complete.
\hfill$\Box$

\section{Strong law of large numbers}

In this section, we give the proofs of Theorems \ref{thrm1} and \ref{main theorem}.
We start with a lemma.

\begin{lemma}\label{lem:U1}
Suppose that Assumption \ref{assum1} holds
 and
$f\in L^2(E;m)\cap L^4(E;m)$ with $\langle f,{\widehat \phi}_0\rangle_m=0$.
Then for any $\mu\in\mathcal{M}_F(E)$ and
$\widetilde{a}\in (0, (\lambda_0-\Re(\lambda_1))\wedge(\lambda_0/2))$,
\begin{equation}\label{5.19}
  \sup_{n>10t_0}e^{(-\lambda_0+\widetilde{a})n}\P_{\mu}\left(\sup_{n\le t\le n+1}|
  \langle U_qf,X_t\rangle|\right)<\infty.
\end{equation}
\end{lemma}

\noindent \textbf{Proof:}
In this proof, we always assume that $n>10t_0$ and $c$ is a positive constant whose value does not depend on $n$ and  may change from one appearance to another.
Define $J_1^f(t):=\langle T_tU_qf,\mu\rangle$ and
$
J_2^f(t):=\int_{t}^\infty e^{-qu}M^{(u)}_t\,du.
$
By \eqref{5.8}, for any $t>0$,
$$
\P_{\mu}\left(\sup_{n\le t\le n+1}|\langle U_qf,X_t\rangle|\right)
\le \sup_{n\le t\le n+1}|J_1^f(t)|+e^{q(n+1)}\P_{\mu}
\left(\sup_{n\le t\le n+1}|J_{2}^f(t)|\right).
$$
First we consider $J_1^f(t)$.
Since $\langle U_qf,{\widehat \phi}_0\rangle=0,$ by \eqref{moment1},
we have $|T_tU_qf|(x)\le ce^{(\lambda_0-\widetilde{a})t}b_{t_0}(x)^{1/2}$.
Thus for $n>10t_0$
\begin{eqnarray}\label{5.18}
 \sup_{n\le t\le n+1}|J_1^f(t)|&\le&  \sup_{n\le t\le n+1}\langle |T_{t}U_qf|,\mu\rangle\nonumber\\
 &\le& c\sup_{n\le t\le n+1}e^{(\lambda_0-\widetilde{a})t}\langle b_{t_0}^{1/2},\mu\rangle\nonumber\\
  &\le& ce^{(\lambda_0-\widetilde{a})n}.
 \end{eqnarray}
Next we deal with $J_{2}^f(t)$.
For $t\in [n,n+1]$,
$$
  J_{2}^f(t)=\int_{t}^\infty e^{-qu}M^{(u)}_{t}\,du=\int_{n}^\infty
  e^{-qu}M^{(u)}_{t}\textbf{1}_{ t<u}\,du.
$$
Thus for $n>10t_0$,
\begin{eqnarray*}
  &&\P_{\mu}\left(\sup_{n\le t\le n+1}|J_{2}^f(t)|\right)
 \le \int_n^{\infty}e^{-qu}
  \P_{\mu}\Big(\sup_{n\le t\le n+1}\Big(\left|M^{(u)}_{t}\right|
  1_{ t<u}\Big)\Big)\,du\nonumber\\
  &\le& 2\int_n^{\infty}e^{-qu}
  \sqrt{\P_{\mu}(M^{(u)}_u)^2}\,du\le
   2\sqrt{c\langle b_{t_0}^{1/2},\mu\rangle}\int_n^\infty
  e^{-qu}e^{(\lambda_0-\widetilde{a})u}\,du\\
  &\le&c(q-\lambda_0+\widetilde{a})^{-1}e^{-(q-\lambda_0+\widetilde{a})n},
\end{eqnarray*}
where the third equality follows from
\eqref{5.22}, \eqref{6.10} and \eqref{moment2}.
It follows that for $n>10t_0$,
\begin{equation}\label{5.7}
  e^{q(n+1)}\P_{\mu}\left(\sup_{n\le t\le  n+1}J_{2}^f(t)\right)\le c
  e^{(\lambda_0-\widetilde{a})n}.
\end{equation}
Combining \eqref{5.18} and \eqref{5.7},
this yields \eqref{5.19}.
The proof is now complete.
\hfill$\Box$

\bigskip

\noindent\textbf{Proof of Theorem \ref{thrm1}:}
Put $\tilde{f}=f-\langle f,{\widehat \phi}_0\rangle_m\phi_0$.
Note that
$$
U_q\phi_0(x)=\int_0^\infty e^{-qt}T_t\phi_0(x)\,dt=\int_0^\infty e^{-qt}
e^{\lambda_0t}\,dt\phi_0(x)=(q-\lambda_0)^{-1}\phi_0(x)
$$
and
\begin{eqnarray}\label{7.41}
  \langle U_qf,{\widehat \phi}_0\rangle_m &=&\int_0^\infty e^{-qt}\langle T_tf,{\widehat \phi}_0\rangle_m\,dt\nonumber\\
  &=&\int_0^\infty e^{-qt}e^{\lambda_0t}\,dt\langle f,{\widehat \phi}_0\rangle_m=(q-\lambda_0)^{-1}\langle f,{\widehat \phi}_0\rangle_m.
\end{eqnarray}
Thus,
$$
U_qf(x)=\langle f,{\widehat \phi}_0\rangle_m U_q\phi_0(x)+U_q(\tilde{f})(x)=
\langle U_qf,{\widehat \phi}_0\rangle_m\phi_0(x)+U_q(\tilde{f})(x).
$$
Hence, to prove \eqref{7.35}, we only need to show that
\begin{equation}\label{7.37}
   e^{-\lambda_0t}\langle U_q(\tilde{f}), X_t\rangle\to 0,
   \quad  \P_{\mu}\mbox{-a.s.}
\end{equation}
Let $M_n:=\sup_{n\le t\le n+1}e^{-\lambda_0t}\Big|\langle U_q(\tilde{f}), X_t\rangle\Big|$.
By \eqref{5.19}, there is a constant $c>0$ so that
$\P_{\mu}M_n\le ce^{-\widetilde{a}n} $
for every  $n>10t_0$.
We conclude by the Borel-Cantelli lemma that
$M_n\to 0$, as $n\to\infty$, $\P_{\mu}$-a.s.,
from which \eqref{7.37} follows immediately.
The proof is now complete.
\hfill$\Box$

For any $f\ge 0$ and   $q>\max\{K,\lambda_0\}$, define
$$
U^{\phi_0}_qf(x)=\int^\infty_0e^{-qt}T^{\phi_0}_tf(x)dt,\quad x\in E,
$$
where $T^{\phi_0}_t$ is defined in \eqref{h-transf}. It is easy to see that
$\phi_0(x)U^{\phi_0}_qf(x)=U_{q+\lambda_0}(\phi_0f)$.

\begin{prop}\label{prop2}
Suppose that Assumptions \ref{assum1} and \ref{assum2} hold.
For any $0\le f\in C_0(E;\R)$ and $\mu\in \mathcal{M}_F(E)$,
\begin{equation}\label{7.26}
  \lim_{t\to\infty}e^{-\lambda_{0}t}\langle\phi_0f, X_t\rangle=\langle
  f\phi_0,{\widehat \phi}_0\rangle_m W_\infty,
   \quad  \P_{\mu}\mbox{-a.s.}
\end{equation}
\end{prop}

\noindent \textbf{Proof:}
By Theorem \ref{thrm1},
\begin{eqnarray*}
\lim_{t\to\infty}e^{-\lambda_{0}t}\langle\phi_0U^{\phi_0}_qf, X_t\rangle
&=&\lim_{t\to\infty}e^{-\lambda_{0}t}\langle U_{q+\lambda_0}(\phi_0f),X_t\rangle\\
&=& \langle U_{q+\lambda_0}(\phi_0f), \wh{\phi}_0\rangle_m W_\infty,
\quad  \P_{\mu}\mbox{-a.s.}
\end{eqnarray*}
According to \eqref{7.41},
$$
\langle U_{q+\lambda_0}(\phi_0f), {\widehat \phi}_0\rangle_m =\frac{1}{q}\langle \phi_0f, {\widehat \phi}_0\rangle_m.
$$
Therefore, for any $q>\max\{K,\lambda_0\}$,
\begin{equation}\label{limit-pot}
\lim_{t\to\infty}e^{-\lambda_{0}t}\langle\phi_0qU^{\phi_0}_qf,
X_t\rangle=\langle f\phi_0,{\widehat \phi}_0\rangle_m W_\infty,
   \quad  \P_{\mu}\mbox{-a.s.}
\end{equation}
Choose a sequence $q_k>\max\{K,\lambda_0\}$ so that $\lim_{k\to\infty}q_k=\infty$.
Put
\begin{eqnarray*}
\Omega^*:&=&\bigcap_{k\geq 1}
\Big\{\lim_{t\to\infty}e^{-\lambda_{0}t}\langle\phi_0q_kU^{\phi_0}_{q_k}f,
X_t(\omega)\rangle=\langle f\phi_0,{\widehat \phi}_0\rangle_m W_\infty(\omega)\Big\}\\
&&\bigcap\Big\{\lim_{t\to\infty}W_t(\omega)=W_\infty(\omega)\Big\}.
\end{eqnarray*}
Then  $\P_{\mu}(\Omega^*)=1$.
Note that, for any $\omega\in\Omega^*$,
\begin{eqnarray*}
&& \left|e^{-\lambda_{0}t}\langle\phi_0q_kU^{\phi_0}_{q_k}f, X_t(\omega)\rangle
-e^{-\lambda_{0}t}\langle\phi_0f, X_t(\omega)\rangle\right|\\
&\le&
e^{-\lambda_0t}\langle \phi_0|q_kU^{\phi_0}_{q_k}f-f|, X_t(\omega)\rangle\\
&\le& \|q_kU^{\phi_0}_{q_k}f-f\|_\infty e^{-\lambda_0t}\langle \phi_0, X_t(\omega)\rangle,
\end{eqnarray*}
where $\|\cdot\|_\infty$ is the $L^\infty$ norm.
Letting $t\to\infty$, we obtain that,
\begin{equation}\label{differ}
\limsup_{t\to\infty}\left|e^{-\lambda_{0}t}\langle\phi_0q_kU^{\phi_0}_{q_k}f,
X_t(\omega)\rangle-e^{-\lambda_{0}t}\langle\phi_0f, X_t\rangle\right|\le
\|q_kU^{\phi_0}_{q_k}f-f\|_\infty W_\infty(\omega).
\end{equation}
By Assumption \ref{assum2},  $\lim_{k\to\infty}\|q_kU^{\phi_0}_{q_k}f-f\|_\infty=0$.
Thus \eqref{differ} implies that, for $\omega\in\Omega^*$,
\begin{equation}\label{limit-differ}
\lim_{k\to\infty}\limsup_{t\to\infty}\left|e^{-\lambda_{0}t}
\langle\phi_0q_kU^{\phi_0}_{q_k}f, X_t(\omega)\rangle-e^{-\lambda_{0}t}
\langle\phi_0f, X_t(\omega)\rangle\right|=0.
\end{equation}
Now, combining \eqref{limit-pot} and \eqref{limit-differ}, we get \eqref{7.26}.
\hfill$\Box$

\bigskip

\noindent\textbf{Proof of Theorem \ref{main theorem}:}
Note that $E_\partial$ is a compact separable metric space.
According to \cite[Exercise 9.1.16(iii)]{Stroock},
$C_b(E_\partial; \R)$,
the space of bounded continuous $\R$-valued functions $f$ on $E$, is separable.
Therefore $C_0(E;\R)$ is also a  separable space.
 Let $\{f_n, n\ge 1\}$ be a
countable dense subset of $C_0(E;\R)$.
Define
\begin{eqnarray*}
 \Omega_0 &:=& \bigcap_{k\ge1}\left\{\omega \in \Omega: \,
\lim_{t\rightarrow\infty}e^{-\lambda_0 t}\langle f_k \phi_0,
X_t\rangle(\omega) = W_\infty(\omega)\int_{E}f_k(y)\phi_0(y){\widehat \phi}_0(y)m(dy)\right\}\\
&&\bigcap \left\{\omega\in\Omega:\lim_{t\to\infty}W_t(\omega)=W_\infty(\omega)\right\}.
\end{eqnarray*}
By Proposition \ref{prop2}, $\mathbb P_\mu(\Omega_0)=1$ for any $\mu\in \mathcal M_F(E)$.

We first consider \eqref{main result2}
on $\{W_\infty>0\}$.
For each $\omega\in\Omega_0\cap \{W_\infty>0\}$ and $t\ge 0$,
we define two probability measures $\nu_t$ and $\nu$ on $D$,
respectively by
$$
\nu_t(F)(\omega)=\frac{e^{-\lambda_0t}\langle 1_F\phi_0,
X_t\rangle(\omega)}
{W_t(\omega)},\, \mbox{ and }\,
\nu(F)=\int_F \phi_0(y){\widehat \phi}_0(y)m(dy),\quad F\in{\cal B}(E).
$$
Note that the measure $\nu_t$ is well-defined for every $t\ge 0$, and  $\nu_t$  and $\nu$ are  probability measures.
By the definition of $\Omega_0$ we know that $\nu_t$ converges weakly to
$\nu$ as $t\to\infty$.
Since
 $\phi_0$ is strictly positive and continuous on $E$,
if
$f$ is a function on $E$ such that $|f|\le c\phi_0$ for some
$c>0$ and that the discontinuity set of $f$ has zero $m$-measure
(equivalently zero $\nu$-measure),
then $g:= f/\phi$ is a bounded
function with the same set of discontinuity. We
thus have
$$
\lim_{t\to\infty}\int_Eg(x)\nu_t({\rm d}x)=\int_Eg(x)\nu({\rm d}x),
$$
which is equivalent to
\begin{eqnarray*}
\lim_{t\rightarrow\infty}e^{-\lambda_0 t}\langle f, X_t\rangle(\omega)
 =W_\infty(\omega)\int_E{\widehat \phi}_0(y)f(y)m({\rm d}y) \quad
\hbox{for } \omega\in\Omega_0\cap \{M_\infty(\phi)>0\}.
\end{eqnarray*}
If $|f|\le c\phi_0$ for some positive constant $c>0$,
\eqref{main result2} holds automatically on
$\{W_\infty=0\}$.
This completes the proof of the theorem.
\hfill$\Box$

\section{Examples}\label{s:example}

In this section we give some examples.
The main purpose is to illustrate the diverse situations where the
main result of this paper can be applied.
We will not try to give the most general examples possible.

\begin{example}
[Super inward Ornstein-Uhlenbeck processes]\label{Inward}
{\rm Let $d\ge 1$, $E=\R^d$. Suppose the spatial motion
$\xi=\{\xi_t,\Pi_{x}\}$ is an
 Ornstein-Uhlenbeck (OU) process on $\R^d$ with infinitesimal generator
$$
\mathcal{L}=\frac{1}{2}\sigma^{2}\Delta-cx\cdot\nabla\mbox{  on }\mathbb{R}^{d},
$$
where $\sigma,\ c>0$.
Without loss of generality, we assume $\sigma=1$. Let
$\varphi(x):=\left(c/\pi \right)^{d/2}e^{-c\|x\|^{2}}$, and
$m(dx)=\varphi(x)dx$. Then $\xi$ is symmetric with respect to
the probability measure $m(dx)$.
Suppose that the branching rate function $\beta(x)=\beta$ is a positive constant,
and the branching mechanism $\psi$ is given by
\begin{equation}\label{bm}
\psi(x,\lambda)=-\lambda+b(x)\lambda^2+\int_{(0,+\infty)}(e^{-\lambda y}-1+\lambda y)n(x,dy),
\quad x\in\R^d, \,\lambda> 0,
\end{equation}
where $b\in \mathcal{B}_b^+(\R^d)$ and $n$ is a kernel from $\R^d$ to $(0,\infty)$ satisfying
$$
  \sup_{x\in \R^d}\int_0^\infty y^2 n(x,dy)<\infty.
$$
Then for the corresponding superprocess,
$$
T_{t}f(x)=e^{\beta t}\Pi_{x}\left[f(\xi_{t})\right]=e^{\beta t}P_tf(x).
$$
It is easy to see that $\lambda_{0}=\beta$, $\phi_0={\widehat \phi}_0=1$ and then $T_t^{\phi_0}=P_t$.

It is well known that, for any $x\in \mathbb{R}^d$, under $\Pi_x$,
$\xi_t$ is of Gaussian distribution with mean $xe^{-ct}$ and variance $\sigma_t^2$,
where $\sigma_t^2:=(1-e^{-2ct})/(2c)$.
The transition density of $\xi_t$ with respect to the probability measure $m(dx)$ on $\R^d$ is given by
$$
p(t,x,y):=
\left(\frac{1}{2c \sigma_t^2}\right)^{d/2}\exp \left(c \| y\|^2 -\frac{\|y-xe^{-ct}\|^2}{2\sigma_t^2} \right).
$$
Note that
$p(t,x, x)=(2\pi \sigma_t^2)^{-d/2}\exp\left(-c \frac{1-e^{-ct}}{1+e^{-ct}} \|x\|^2\right) /\varphi (x)$.
Thus $a(t)= p(2t, x, x)$ is $L^1(\R^d; m)$-integrable for all $t>0$
and there is some $t_0>0$ so that $a(t)\in L^2(\R^d; m)$-integrable for $t\geq t_0$.
Hence Assumption \ref{assum1} holds for $\xi$.

For any $f\in C_0(\R^d;\R)$, we have
$$
P_tf(x)=\int_{\R^d}p_t(x,y)f(y) m(dy)
=\int_{\R^d}(2\pi)^{-d/2}\exp\left(-\|y\|^2/2\right)
f(\sigma_ty+xe^{-c t})dy.
$$
Using the dominated convergence theorem, one can easily check that
$P_tf\in C_0(\R^d;\R)$.
Suppose $f$ is a continuous function with compact support.
Let  $M_0>0$ so that $f(x)=0$ for $\|x\|\geq M_0$.
For any $M>0$,
 $$\begin{array}{rl}\left|P_tf(x)-f(x)\right|
 =&\displaystyle\left|\int_{\R^d}(2\pi)^{-d/2}\exp\left(-\|y\|^2/2\right)
 \left[f(\sigma_ty+xe^{-c t})-f(x)\right]dy\right|\\
 \le&\displaystyle\int_{\R^d}(2\pi)^{-d/2}\exp\left(-\|y\|^2/2\right)
 \left|f(\sigma_ty+xe^{-c t})-f(x)\right|dy\\
 \le&\displaystyle\int_{\|y\|\le M}(2\pi)^{-d/2}\exp\left(-\|y\|^2/2\right)
 \left|f(\sigma_ty+xe^{-c t})-f(x)\right|dy\\
 &\displaystyle+2\|f\|_\infty\int_{\|y\|\ge M}(2\pi)^{-d/2}\exp\left(-\|y\|^2/2\right)\,dy\\
 =:&I+II.\end{array}$$
For any $\epsilon>0$, we choose $M>0$ such that $II\le \epsilon/2$.
For part $I$,
we claim that, for any $\epsilon>0$, there exists $\delta$, for $t\le \delta$,$$\sup_{\|y\|\le M}\sup_{x\in \R^d}\left|f(\sigma_ty+xe^{-c t})-f(x)\right|\le \epsilon/2.$$
Therefore $I<\epsilon/2$, and then $\|P_tf-f\|_\infty\to 0$ as $t\to 0$.

Now we prove the claim. Note that $$\left|f(\sigma_ty+xe^{-c t})-f(x)\right|\le \left|f(\sigma_ty+xe^{-c t})-f(xe^{-ct})\right|+\left|f(xe^{-ct})-f(x)\right|.$$
Since $f$ is  uniformly continuous on $\R^d$,
there is a constant $\delta_0>0$ such that
$|f(y)-f(x)|\le \epsilon/4$
for any $x,y$ satisfying $\|x-y\|\le \delta_0$.
Since $\sigma_t\to 0$ as $t\to 0$, there exists $\delta_1>0$ such that, for $t<\delta_1$,
$\|\sigma_t\|\le \delta_0/M$, and then $\sup_{\|y\|\le M}\sup_{x\in \R^d}|f(\sigma_ty+xe^{-c t})-f(xe^{-ct})|\le \epsilon/4$.
Choose $\delta_2$, such that for $t\le \delta_2$, $e^{ct}-1\le \delta_0/M_0.$
Then , for $t\le \delta_2$, $$|f(xe^{-ct})-f(x)|\le |f(xe^{-ct})-f(x)|\textbf{1}_{\|x\|\le M_0e^{ct}}\le \epsilon/4,$$
where in the second inequality we use the fact that $\|xe^{-ct}-x\|=\|x\|(1-e^{-ct})\le M_0(e^{ct}-1)\le \delta_0$.
Then, choosing $\delta=\delta_1\wedge\delta_2$, we prove the claim.

For general $f\in C_0(\R^d;\R)$,
there exist continuous functions $f_n$ with compact support such that $\|f_n-f\|_\infty\to0$, as $n\to\infty$.
Then
\begin{eqnarray*}\|P_tf-f\|_\infty&\le& \|P_tf-P_tf_n\|_\infty+\|P_tf_n-f_n\|_\infty+\|f_n-f\|_\infty\\
&\le&\|P_tf_n-f_n\|_\infty+2\|f_n-f\|_\infty.
\end{eqnarray*}
Letting $t\to0$ and then $n\to\infty$, we get that  $\|P_tf-f\|_\infty\to 0$ as $t\to 0$.
Since $T^{\phi_0}_t=P_t$,  Assumption \ref{assum2} is satisfied. Therefore
for the superprocess in this example, all our assumptions are satisfied.

This example covers Examples 4.1 and 4.6 in \cite{EKW14}.
For variable $\alpha (x)=\beta (x) a(x)$, see Example \ref{E:4.9}.
}
\end{example}

\begin{example}\label{E:4.2}
[Super outward Ornstein-Uhlenbeck processes]\label{Outward}
{\rm Let $d\ge 1$, $E=\R^d$.
Suppose the spatial motion $\xi=\{\xi_t,\Pi_{x}\}$ is an
OU process on $\R^d$ with infinitesimal generator
$$
\mathcal{L}=\frac{1}{2}\sigma^{2}\Delta+cx\cdot\nabla\mbox{  on }\mathbb{R}^{d},
$$
where $\sigma,\ c>0$.
Without loss of generality, we assume $\sigma=1$.
Under $\Pi_x$,
$\xi_t$ is of Gaussian distribution with mean $xe^{ct}$ and variance $(e^{2ct}-1)/(2c)$.

Let $\tilde{\varphi}(x):=\left(c/\pi\right)^{-d/2}e^{c\|x\|^{2}}$, and
$m(dx)=\tilde{\varphi}(x)dx$. Then $\xi$ is symmetric with respect to
the $\sigma$-finite measure $m(dx)$.
As in the previous example, we suppose that the branching rate function $\beta(x)=\beta$
is a positive constant, and the branching mechanism $\psi$ is given by
\eqref{bm}.
Then for the corresponding superprocess,
$$
T_{t}f(x)=e^{\beta t}\Pi_{x}\left[f(\xi_{t})\right]=e^{\beta t}P_tf(x).
$$
The generator of $\{T_t: t\ge 0\}$ is $\mathcal{L}+\beta$.

The transition density of $\xi$ with respect to the measure $m$ is
$$
p(t,x,y)=
\left(\frac{1}{e^{2ct}-1}\right)^{d/2}\exp\left( -\frac{c}{(1-e^{-2ct})}
\left(\|y\|^2+\|x\|^2-2x\cdot ye^{-ct}\right)\right).
$$
Thus
$$
a_t(x)=p(2t,x,x)=\left(\frac{1}{e^{2ct}-1}\right)^{d/2}\exp
\left(-\frac{2c\|x\|^2}{(1+e^{-ct})} \right).
$$
It is obvious that $a_t\in L^1(\R^d;m)\cap L^2(\R^d;m)$. Thus Assumption \ref{assum1} is satisfied.
Suppose $\beta(x)=\beta\in(cd,\infty)$.

The operator ${\cal L}+cd$ is the formal adjoint of the
inward OU process
with infinitesimal generator  $\frac{1}{2}\sigma^{2}\Delta-cx\cdot\nabla$
on $\mathbb{R}^{d}.$ Since $\varphi(x)$ defined in Example \ref{Inward} is
the invariant density of $\frac{1}{2}\sigma^{2}\Delta-cx\cdot\nabla$
on $\mathbb{R}^{d}$, $({\cal L}+cd)\varphi=0$. Thus we have $({\cal L}+\beta)\varphi=(\beta-cd)\varphi$.
Since $\varphi\in L^2(\R^d,m)$ and $\varphi$ is strictly positive everywhere,
we know that $\phi_0={\widehat \phi}_0=\varphi$ and $\lambda_0=\beta-cd.$ Thus
$$T_t^{\phi_0}f(x)=\frac{e^{cdt}P_t(f\varphi)(x)}{\varphi(x)}=\tilde{P}_tf(x),$$
where $\tilde{P}_t$  is the semigroup of the inward OU-process with infinitesimal generator
$$
\frac{1}{2}\Delta-cx\cdot\nabla\mbox{  on }\mathbb{R}^{d}.
$$
From the discussion  in Example \ref{Inward}, we see that Assumption \ref{assum2} is satisfied.
Thus, when $\beta(x)=\beta\in(cd,\infty)$, the superprocess of this example satisfies
all our assumptions.

This example covers Examples 4.2 in \cite{EKW14}.}
\end{example}

\begin{example}{\rm
 Suppose that $\eta=\{\eta_t,\Pi_x\}$ is an $m$-symmetric Hunt process on $E$
 and that $\eta$ has a transition density $\widetilde{p}(t,x,y)$ with respect to $m$.
 Suppose also that $\widetilde{p}$ is strictly positive, continuous and satisfies Assumption \ref{assum1}.
 Let $\{\widetilde P_t, t\geq 0\}$ be the transition semigroup of $\eta$ on $L^2(E; m)$.
 Since, for each $t>0$, $\wt P_t$ is compact,
 the infinitesimal generator $\wt \sL$ of $\{\wt P_t, t\geq 0\}$ has
 discrete spectrum: $0\geq \wt \lambda_0>\wt \lambda_1 \geq \cdots $.
 Denote   the corresponding normalized eigenfunctions by $\{\wt \phi_k; k\geq 0\}$, with
 $\| \wt \phi_k \|_{L^2(E; m)}=1$ for every $k\geq 0$.
 We can choose $\wt \phi_0$ so that it is strictly positive and continuous.
 By the spectral representation, we can express $\wt p(t, x, y)$ by $\sum_{k=0}^\infty e^{\wt \lambda_k t} \wt \phi_k (x) \wt \phi_k (y)$.
 It follows that $\wt p(t, x, x)$ is decreasing in $t>0$; see \cite[Section 2]{DS}.
  Define
 $$
 \widetilde{P}_t^{\widetilde{\phi}_0}f:=e^{-\widetilde{\lambda}_0t}\frac{\widetilde{P}_t(f\widetilde{\phi}_0)(x)}
 {\widetilde{\phi}_0(x)}.
 $$
Assume
that $\widetilde{P}_t^{\widetilde{\phi}_0}$ satisfies  Assumption \ref{assum2}.

Let $S_t$ be a subordinator, independent of $\xi$, with drift $b>0$. Then $S_t\ge bt$. Let $\phi$  be the Laplace exponent of $S$, that is,
$$
 \E(e^{-\theta S_t})=e^{-t\phi(\theta)}, \qquad \theta>0.
$$
Suppose that $\alpha(x)=\alpha$ is a constant function
and satisfies $\alpha>\phi(-\widetilde{\lambda}_0)$.
We put $\xi_t:=\eta_{S_t}$. Let $P_t$ be the semigroup of $\xi$ and $p(t,x,y)$ be the transition density of $\xi$ with respect to $m$.
Then $p(t,x,y)=\E \widetilde{p}(S_t,x,y)$.
Since $t\to \widetilde{p}(t,x,x)$ is a decreasing function,
 $p(2t,x,x)=\E \widetilde{p}(S_{2t},x,x)\le \widetilde{p}(2bt,x,x)$,
which implies that $\eta$ satisfies Assumption \ref{assum1}.
Note that $T_t=e^{\alpha t}P_t$, and
$$P_t\widetilde{\phi}_0(x)=\E(\widetilde{P}_{S_t}\widetilde{\phi}_0(x))
=\E  e^{\widetilde{\lambda}_0S_t}\widetilde{\phi}_0(x)=e^{-t\phi(-\widetilde{\lambda}_0)}\widetilde{\phi}_0(x).$$
Thus, $\lambda_0=\alpha-\phi(-\widetilde{\lambda}_0)>0$ and $\phi_0=\widetilde{\phi}_0$.
Then $$T_t^{\phi_0}f(x)=e^{t\phi(-\widetilde{\lambda}_0)}\frac{P_t(f\phi_0)(x)}{\phi_0(x)}
=e^{t\phi(-\widetilde{\lambda}_0)}\E\left[\frac{\widetilde{P}_{S_t}(f\phi_0)(x)}{\phi_0(x)}\right].$$
Thus, we have
\begin{eqnarray*}|T_t^{\phi_0}f(x)-f(x)|&\le& e^{t\phi(-\widetilde{\lambda}_0)}\E\left|\frac{\widetilde{P}_{S_t}(f\phi_0)(x)}{\phi_0(x)}-e^{\widetilde{\lambda}_0S_t}f(x)\right|\\
&\le& e^{bt\widetilde{\lambda}_0}e^{t\phi(-\widetilde{\lambda}_0)}
\E\left[\|\widetilde{P}_{S_t}^{\widetilde{\phi}_0}f-f\|_\infty\right].
\end{eqnarray*}
Since $\|\widetilde{P}_{S_t}^{\widetilde{\phi}_0}f-f\|_\infty\to 0,$ as $t\to0$, and $\|\widetilde{P}_{S_t}^{\widetilde{\phi}_0}f-f\|_\infty\le 2\|f\|_\infty$, using the dominated convergence theorem, we get that
$$\lim_{t\to0}\|T_t^{\phi_0}f-f\|_\infty=0.$$
Thus, the superprocess of this example satisfies
all our assumptions.

In particular, this example is applicable when $\eta$ is the outward
Ornstein-Uhlenbeck process or inward  Ornstein-Uhlenbeck process dealt with in the
Examples \ref{Inward} and \ref{Outward}.
}
\end{example}

The next two examples give the cases when $\alpha$ is not a constant function.

\begin{example}[Pure jump SBM]\label{S:2.1} {\rm
Suppose that $S=\{S_t, t\ge 0\}$ is a drift-free subordinator.
The Laplace exponent $\phi$ of $S$ can be written in the form
\begin{equation}\label{e:LK}
\phi(\lambda)=\int^\infty_0(1-e^{-\lambda t})\,u(dt),
\end{equation}
where $u$ is a measure on $(0, \infty)$ satisfying $\int^\infty_0
(1\wedge t)\,u(dt)<\infty$.
The measure $u$ is
the L\'evy measure of the subordinator (or of $\phi$).
In this example, we will assume that $\phi$ is a complete Bernstein function, that
is, the measure $u$ has a completely monotone density, which we also denote by $u$.

Let $W=\{W_t, t\ge 0\}$
 be a Brownian motion in $\R^d$ independent of the subordinator $S$.
The subordinate Brownian motion
$Y=\{Y_t, t\ge 0\}$ is defined by
$Y_t:=W_{S_t}$, which is a rotationally symmetric L\'evy process with L\'evy exponent
$\phi( |\xi|^2)$.
It is known that the L\'evy measure of the process $Y$ has a density given by
$x \to j(|x|)$ where
\begin{equation}\label{e:ld4s}
j(r):=\int^{\infty}_0(4\pi t)^{-d/2}e^{-r^2/(4t)}\,u(t)dt, \qquad r>0.
\end{equation}
 Note that the function
$r\mapsto j(r)$ is continuous and decreasing
on $(0, \infty)$.

Suppose that $\phi$ satisfies the following
growth condition at infinity:

\medskip
\noindent
{\bf (A):}
There exist constants $ \delta_1, \delta_2 \in (0,1)$, $a_1\in (0, 1)$, $a_2\in (1, \infty)$
and $R_0>0$ such that
\begin{eqnarray*}
a_1\lambda^{\delta_1} \phi(r) \leq
\phi(\lambda r) \le a_2 \lambda^{\delta_2} \phi(r) \quad
&\hbox{for } \lambda \ge 1 \hbox{ and }  r \ge R_0.
\end{eqnarray*}
See \cite{CKS1} for examples of a large class of symmetric L\'evy processes
satisfying condition {\bf (A)}.

Suppose $D$ is a bounded $C^{1,1}$ open set with characteristics
$(R_0, \Lambda)$, and let $\xi$ be the subprocess of $Y$ killed upon
leaving $D$. It is known that $\xi$ is a Feller process with strong Feller property in $D$.
Moveover, by \cite[Corollary 1.6]{CKS1}, $\xi$ has a jointly continuous transition density
function $p_D(t, x, y)$ with respect to the Lebesgue measure on $D$ so that
for every $T>0$, there exist
$c_1=c_1(R_0, \Lambda, T, d, \phi)\geq 1$ and $c_2=c_2(R_0, \Lambda, T, d, \phi)>0$
 such that for $0<t\leq T$, $x,y\in D$,
 \begin{eqnarray}
     &&  c_1^{-1} \left(1\wedge \frac{\Phi (\delta_D(x))}t \right)^{1/2} \left(1\wedge
 \frac{\Phi (\delta_D(y))}t \right)^{1/2}
 \left( \Phi^{-1}(t)^{-d} \wedge t j(|x- y|)\right)\nonumber\\
 & \le & {p_{D}(t, x, y)}\label{e:4}\\
 & \le & c_2 \left(1\wedge \frac{\Phi (\delta_D(x))}t \right)^{1/2} \left(1\wedge
 \frac{\Phi (\delta_D(y))}t \right)^{1/2}\left( \Phi^{-1}(t)^{-d} \wedge t j(c_2|x- y|/4)\right).\nonumber
\end{eqnarray}
Here $\Phi(r):=\frac1{\phi(r^{-2})}$,
$j$ is the function defined in \eqref{e:ld4s},
and $\delta_D(x)$ is the Euclidean distance between $x$ and $\partial D$.
Since $p_D(t,x,y)$ is symmetric, $a_t(x)=p_D(2t,x,x)\le c_2\Phi^{-1}(2t)^{-d}$.
Thus, Assumption \ref{assum1} is satisfied.

Suppose that the branching rate function $\beta$ and the branching mechanism
satisfy the assumptions of Subsection \ref{superp}, and that
the corresponding superprocess $X$ is supercritical. The corresponding semigroup
$\{T_t:t\ge 0\}$ has a continuous density $q(t, x, y)$
satisfying the same two-sided estimates \eqref{e:4} with possibly different $c_1\geq 1$ and $c_2$.
Since $\phi_0(x)=e^{\lambda_0t}T_t\phi_0(x)$, by \eqref{e:4},
$$
\phi_0(x) \asymp \Phi (\delta_D (x))^{1/2}.
$$
We now show that Assumption \ref{assum2} holds.
Suppose $f\in C_0 (D)$.
For any given $\eps >0$, there $\delta>0$ so that
$|f(x)-f(y)|< \eps$ whenever $|x-y|<\delta$. Hence by the display above and
\eqref{e:4}, for small $t>0$,
\begin{eqnarray}
&&\sup_{x\in D} |T^{\phi_0}_t f(x)-f(x)|\nonumber\\
&=& \sup_{x\in D}  \frac{e^{-\lambda_0t} \left| \Pi_x \left[ e^{\int_0^t\alpha(\xi_s)\,du}  \phi_0(\xi_t)
\left( f(\xi_t)- f(\xi_0)\right)\right] \right|}{\phi_0(x)}  \nonumber\\
&\leq & \eps +\sup_{x\in D}  \frac{e^{-\lambda_0t} \left| \Pi_x \left[ e^{\int_0^t\alpha(\xi_s)\,du}  \phi_0(\xi_t)
\left| f(\xi_t)- f(\xi_0)\right|; |\xi_t-\xi_0|\geq  \delta\right] \right|}{\phi_0(x)}
\nonumber\\
&\leq & \eps + \sup_{x\in D}  c \frac{e^{(-\lambda_0+\| \alpha\|_\infty) t} \|\phi_0\|_\infty
\|f\|_\infty  \Pi_x ( |\xi_t-\xi_0| \geq  \delta) }{\phi_0(x)} \nonumber\\
&\leq & \eps + \sup_{x\in D}   c \frac{   \Phi (\delta_D (x))^{1/2} t^{-1/2}
  \int_{y\in D: |y-x|>\delta} t j(c_2|y-x|/4) dy} { \Phi (\delta_D (x))^{1/2}}
  \nonumber \\
&\leq & \eps + c\sqrt{t}
\int_{|z|\geq c_2\delta/4} (1\wedge |z|^2) j(|z|) dz.
\end{eqnarray}
It follows that $\lim_{t\to 0} \| T^{\phi_0}_t f -f \|_\infty=0$ and
Assumption \ref{assum2} is satisfied.
}
\end{example}

\begin{example}[SBM with Gaussian component]\label{S:2.2}
{\rm
Suppose that $S=\{S_t, t\ge 0\}$
is a subordinator with drift $b>0$.
The Laplace exponent $\phi$ of $S$ can be written in the form
\begin{equation}\label{e:LK2}
\phi(\lambda)=b\lambda+\int^\infty_0(1-e^{-\lambda t})\,u(dt),
\end{equation}
where $u$ is a measure on $(0, \infty)$ satisfying $\int^\infty_0
(1\wedge t)\,u(dt)<\infty$.
Without loss of generality we assume that $b=1$.
In this example, we will assume that $\phi$ is a complete Bernstein function
and that the L\'evy density $u(t)$ of $S$ satisfies
the following growth condition on
$u(t)$ in \eqref{e:LK} near zero:
For any $M>0$, there exists $c=c(M)>1$ such that
\begin{equation}\label{H:1a}
u(r)\le c u(2r), \quad r\in (0, M).
\end{equation}
Let $W=\{W_t, t\ge 0\}$
 be a Brownian motion in $\R^d$ independent of the subordinator $S$.
The subordinate Brownian motion
$Y=\{Y_t, t\ge 0\}$ is defined by
$Y_t:=W_{S_t}$,
which is a rotationally symmetric L\'evy process with L\'evy exponent $\phi
( |\xi|^2)$.
It is known that the L\'evy measure of the process $Y$
has a density $j(|x|)$ given by \eqref{e:ld4s}.

For any open set $D\subset \R^d$ and positive constants $c_1$ and $c_2$, we define
\begin{eqnarray}\label{eq:qd}
&&h_{D, c_1, c_2}(t, x, y)\\
&:=&
 \left(1\wedge \frac{\delta_D(x)}{\sqrt{t}}\right)
\left(1\wedge \frac{\delta_D(y)}{\sqrt{t}}\right)
\left(t^{-d/2}e^{-c_1|x-y|^2/t}+ t^{-d/2}\wedge (tj(c_2|x-y|)\right).\nonumber
\end{eqnarray}

Suppose $D$ is a bounded $C^{1,1}$ open set with characteristics
$(R_0, \Lambda)$, and let $\xi$ be the subprocess of $Y$ killed upon
leaving $D$.
It is known that $\xi$ is a Hunt process symmetric with respect to the
Lebesgue measure on $D$ and that $\xi$ has a strictly positive
continuous transition density
$p_D(t, x, y)$ with respect to the Lebesgue
measure on $D$.
We assume the following upper bound condition on the transition
density function $\widetilde{p}(t, |x|)$ of $Y$:
for any $T>0$, there exist $C_j
\ge 1$, $j=1, 2, 3, $ such that
for all $(t,r)\in (0, T]\times [0, $ diam$(D)]$,
\begin{equation}\label{globalup}
\widetilde{p}(t, r) \le C_1 \left(t^{-d/2}e^{-r^2/C_2t}+ t^{-d/2}\wedge (tj(r/C_3))\right).
\end{equation}
It is established in \cite{CK3} that the  above estimate holds for
a large class of symmetric diffusion processes with jumps with $D=\R^d$.
Using Meyer's method of removing and adding jumps, it can be shown
that \eqref{globalup} is true for a larger class of symmetric Markov processes,
including subordinate Brownian motions with Gaussian components under some
additional condition. See the paragraph containing (1.12) in \cite{CKS2} for more information.

The following is proved in \cite[Theorem 1]{CKS2}.
\begin{description}
\item{\rm (i)}
For every $T>0$, there exist
  $c_1=c_1(R_0, \Lambda_0,   \lambda_0, T,
  \psi, d)>0$ and  $c_2=c_2(R_0, \Lambda_0,   \lambda_0,  d)>0$
such that for all $(t, x, y) \in (0, T]\times D\times D$,
\begin{equation}\label{e:2.9}
p_D(t, x, y)\,\ge\, c_1 \, h_{D, c_2, 1}(t, x, y).
\end{equation}
\item{\rm (ii)}
If $D$ satisfies \eqref{globalup}, then for every $T>0$,
there exists \\$c_3=
c_3(R_0, \Lambda_0,  T, d,  \psi, C_1, C_2, C_3, d )>1$ such
that for all $(t, x, y) \in (0, T]\times D\times D$,
\begin{equation}\label{e:2.10}
p_D(t, x, y)\,\le\, c_3\,
h_{D, C_4, C_5}(t, x, y),
\end{equation}
where  $C_4=(16C_2)^{-1}$ and
$C_5= (8 \vee 4C_3)^{-1}$.
\end{description}

Let $E=D$ and $m$ be the Lebesgue measure on $D$.
Since $p_D(t,x,y)$ is symmetric, $a_t(x)=p_D(2t,x,x)\le c t^{-d/2}$.
Thus, Assumption \ref{assum1} is satisfied.

Suppose that the branching rate function $\beta$ and the branching mechanism
satisfy the assumptions of Subsection \ref{superp}, and that
the corresponding superprocess $X$ is supercritical.
Using the above two-sided heat kernel estimate for $\xi$, we can establish in a similar way
as in Example  \ref{S:2.1} that Assumption \ref{assum2} also holds.
}
\end{example}

\begin{remark}\rm
In fact, in the two examples above,
$\xi$ does not need to be a subordinate Brownian motion killed upon leaving $D$.
All we need are the heat kernel estimates like \eqref{e:4} or \eqref{e:2.9}-\eqref{e:2.10}.
For example, suppose $Y^D$ is the subprocess of some subordinate Brownian motion $Y$ killed upon leave $D$
 that has the property \eqref{e:4} or \eqref{e:2.9}-\eqref{e:2.10}. Let $\xi$ be a Markov process
obtained from $Y^D$ though a  Feynman-Kac transform with
bounded potential function. Then
$\xi$ enjoys the property  \eqref{e:4} or \eqref{e:2.9}-\eqref{e:2.10}.
For other examples of processes that satisfy two-sided bounds similar to \eqref{e:4}, including censored stable processes in $C^{1,1}$ open sets and their
 local and non-local Feynman-Kac transforms, see \cite{CKS4}.
Our main results are applicable to these processes as well.
\end{remark}

In all the examples above, the spatial motion $\xi$ is symmetric. Now we give two examples
where the spatial motion $\xi$ is not symmetric.

\begin{example}\label{newex1}
{\rm
Suppose $d\ge 3$ and that $\nu=(\nu^1, \cdots, \nu^d)$, where each $\nu^j$ is a signed measure
on
$\mathbb{R}^d$ such that
$$
 \lim_{r\to 0}\sup_{x\in\mathbb{R}^d}\int_{B(x, r)}\frac{|\nu^j|(dy)}{|x-y|^{d-1}}=0.
$$
Let $\xi^{(1)}=\{\xi^{(1)}_t, t\ge 0\}$
be a Brownian motion with drift $\nu$ in $\mathbb{R}^d$,
see \cite{BC}.
Suppose that $D$ is a bounded  domain in $\mathbb{R}^d$.
Let $M>0$ so that $B(0, M/2)\supset D$.
 Put $B=B(0, M)$. Let $G_B$ be the Green function of $\xi^{(1)}$ in $B$ and define
$H(x):=\int_BG_B(y, x)dy$.
Then $H$ is a strictly positive continuous function on $B$.
Let $\xi$ be the
process obtained by killing $\xi^{(1)}$ upon exiting $D$.
$\xi$ is a Hunt process and it has a strictly positive continuous
transition density $\widetilde{p}(t, x, y)$ with respect to
the Lebesgue measure on $D$.
Let $E=D$ and $m$ be the measure defined by $m(dx)=H(x)dx$.
It follows from \cite{KiSo08, KiSo08c} that
$\xi$ has a dual process with respect to $m$.
The transition density of $\xi$ with respect to $m$ is
given by $p(t, x, y)=\widetilde{p}(t, x, y)/H(y)$.

Suppose further that $D$ is $C^{1,1}$, then it follows from
\cite[Theorem 4.6]{KiSo06} that there exist $c_1>1, c_2>c_3>0$ such that
for all $(t, x, y)\in (0, 1]\times D\times D$,
 \begin{eqnarray*}
     &&  c_1^{-1} t^{-d/2}\left(1\wedge \frac{\delta_D(x)}{\sqrt{t}}\right)
     \left(1\wedge \frac{\delta_D(y)}{\sqrt{t}}\right)
     \exp\left(-\frac{c_2|x-y|^2}{t}\right)\\
 & \le & {\tilde{p}(t, x, y)}
    \le c_1t^{-d/2}\left(1\wedge \frac{\delta_D(x)}{\sqrt{t}}\right)
     \left(1\wedge \frac{\delta_D(y)}{\sqrt{t}}\right)
 \exp\left(-\frac{c_3|x-y|^2}{t}\right).
\end{eqnarray*}

It follows from the display above and the semigroup property that,
for any $t>0$, $\tilde{p}(t, x, y)$ is bounded.
By \cite[(2.6)]{KiSo08c}, $H(x)\asymp \delta_B(x)$. So for $x\in D,$ $c\le H(x)\le C$, where $c,C>0$.
Thus, $p(t,x,y)$ is also bounded in $D$ and $m$ is a finite measure.
Thus Assumption \ref{assum1}
is satisfied.

Suppose that the branching rate function $\beta$ and the branching mechanism
satisfy the assumptions of Subsection \ref{superp}, and that
the corresponding superprocess $X$ is supercritical.
Using the above two-sided heat kernel estimate for $\xi$, we can establish in a similar way
as in Example  \ref{S:2.1} that Assumption \ref{assum2} also holds.
}
\end{example}

\begin{example}\label{newex2}
{\rm
Suppose $d\ge 2$, $\alpha\in (1, 2)$,
and that $\nu=(\nu^1, \cdots, \nu^d)$, where each $\nu^j$ is a signed measure
on $\mathbb{R}^d$ such that
$$
 \lim_{r\to 0}\sup_{x\in\mathbb{R}^d}\int_{B(x, r)}\frac{|\nu^j|(dy)}{|x-y|^{d-\alpha+1}}=0.
$$
Let $\xi^{(2)}=\{\xi^{(2)}_t, t\ge 0\}$
be an $\alpha$-stable process with drift
$\nu$ in $\mathbb{R}^d$, see \cite{KiSo13}.
Suppose that $D$ is a bounded open set in $\mathbb{R}^d$ and suppose $M>0$ is such that
$D\subset B(0, M/2)$. Put $B=B(0, M)$.
Let $G_B$ be the Green function of $\xi^{(2)}$ in $B$ and define $H(x):=\int_BG_B(y, x)dy$.
Then $H$ is a strictly positive continuous function on $B$.
Let $\xi$ be the
process obtained by killing $\xi^{(2)}$ upon exiting $D$.
$\xi$ is a Hunt process and it has a strictly positive continuous
transition density $\widetilde{p}(t, x, y)$ with respect to
the Lebesgue measure on $D$.
Let $E=D$ and $m$ be the measure defined by $m(dx)=H(x)dx$.
It follows from \cite[Section 5]{CKS3} and \cite{KiSo13} that
$\xi$ has a dual process with respect to $m$.
The transition density of $\xi$ with respect to $m$ is
given by $p(t, x, y)=\widetilde{p}(t, x, y)/H(y)$.
By \cite[Corollary 1.4]{CKS3} and \cite{KiSo14}, we can check that $H(x)\asymp \delta_B(x)^{\alpha/2}$.
Thus, for $x\in D,$ $c\le H(x)\le C$, for some $c,C>0$.

Suppose further that $D$ is $C^{1,1}$, then it follows from
\cite[Theorem 1.3]{CKS3} and \cite{KiSo14} that there exists $c_1>1$ such that
for all $(t, x, y)\in (0, 1]\times D\times D$,
 \begin{eqnarray*}
     &&  c_1^{-1} \left(1\wedge \frac{\delta^{\alpha/2}_D(x)}{\sqrt{t}}\right)
     \left(1\wedge \frac{\delta^{\alpha/2}_D(y)}{\sqrt{t}}\right)
 \left(t^{-d/\alpha} \wedge \frac{t}{|x-y|^{d+\alpha}}\right)\\
 & \le & {\tilde{p}(t, x, y)}
    \le c_1\left(1\wedge \frac{\delta^{\alpha/2}_D(x)}{\sqrt{t}}\right)
     \left(1\wedge \frac{\delta^{\alpha/2}_D(y)}{\sqrt{t}}\right)
 \left(t^{-d/\alpha} \wedge \frac{t}{|x-y|^{d+\alpha}}\right).
\end{eqnarray*}
It follows from the display above and the semigroup property that,
for any $t>0$, $\tilde{p}(t, x, y)$ is bounded.
Since $H$ is bounded between two positive constants,
Assumption \ref{assum1} is satisfied.

Suppose that the branching rate function $\beta$ and the branching mechanism
satisfy the assumptions of Subsection \ref{superp}, and that
the corresponding superprocess $X$ is supercritical.
Using the above two-sided heat kernel estimate for $\xi$, we can establish in a similar way
as in Example  \ref{S:2.1} that Assumption \ref{assum2} also holds.
}
\end{example}

In the following example, our main result does not apply directly.
However, we could apply our main result after a transform.

\begin{example}\label{E:4.9}
{\rm Suppose the spatial motion
$\xi=\{\xi_t,\Pi_{x}\}$
is an OU-process on $\R^d$ with
infinitesimal generator
$$\mathcal{L}=\frac{1}{2}\sigma^{2}\Delta-cx\cdot\nabla\mbox{  on }\mathbb{R}^{d},$$
where $\sigma,\ c>0$.
Without loss of generality, we assume $\sigma=1$. Let $\varphi(x):=
\left(c/\pi\right)^{d/2}e^{-c|x|^{2}}$, and
$m(dx)=\varphi(x)dx$. Then $(\xi,\Pi_{x})$ is symmetric with respect to
the probability measure $m(dx)$.

Let  $a(x)=c_1|x|^{2}+c_2$ with $c_1,c_2>0$, and let $P^{a}_{t}$ be the Feynman-Kac semigroup,
$$P^{a}_{t}f(x):=\Pi_{x}\left[e^{\int_{0}^{t}a(\xi_{s})ds}f(\xi_{t})\right].$$
Suppose $c>\sqrt{2c_1}$ and write $\upsilon=\frac{1}{2}(c-\sqrt{c^{2}-2c_1})$. Let
$$\lambda_{c}:=\inf\{\lambda\in\mathbb{R}:\ \mbox{there exists }u>0
\mbox{ such that }(\mathcal{L}+a-\lambda)u=0\mbox{ in }\mathbb{R}^{d}\}$$
be the generalized principal eigenvalue. Let $h$ denote the corresponding
ground state, i.e., $h>0$ such that
$(\mathcal{L}+a-\lambda_{c})h=0$. As is indicated
in \cite{Englander2}, $\lambda_{c}=c_2+d\upsilon>0$ and
$h(x)=\left(\frac{c-2\upsilon}{c}\right)^{d/2}\exp\{\upsilon\|x\|^{2}\}$.

Note that $h=e^{-\lambda_{c}t}P^{a}_{t}h$ on $\mathbb{R}^{d}$.
Let $\Pi^{h}_{x}$ be defined as in \eqref{h-transf} with $\phi_0$
replaced by $h$. The transformed process $(\xi, \Pi^{h}_{x})$ is also an
OU-process with infinitesimal generator
$\frac{1}{2}\Delta-(c-2v) x\cdot\nabla$   on $\mathbb{R}^{d}.$

Let $\psi(x,z)=-a(x)z+\alpha(x)z^2$, where $\alpha\in C^{\eta}(\R^d)$,
$\alpha(x)>0$ for all $x\in\R^d$.
A superprocess $X$ with spacial motion $\xi$, branching rate $\beta(x)=1$ and
branching mechanism $\psi$ can be defined by $X=\frac{1}{h}X^h$, where $X^h$
is the superprocess with spacial motion $(\xi, \Pi^{h}_{x})$, branching rate
$\beta(x)=1$ and branching mechanism $\psi^h(x,z)=-\lambda_cz+h(x)\alpha(x)z^2.$

Assume that $h\alpha$ is bounded in $\R^d$. Then, for $X^h$, we have  $m^h(dx)=(\frac{c-2\upsilon}{\pi})^{d/2}e^{-(c-2\upsilon)|x|^2}\,dx$,
$\lambda^h_0=\lambda_c$ and $\phi^h_0=1$.
From the discussion  in Example \ref{Inward}, we see that the Assumption
\ref{assum1}  and Assumption \ref{assum2} are satisfied  for the superprocess $X^h$.
Then, there
exists $\Omega_0\subset\Omega$ of probability one (that is,
$\mathbb P_\mu (\Omega_0)=1$ for every $\mu\in \mathcal M_F(\R^d)$) such that, for
every $\omega\in\Omega_0$ and for every bounded Borel measurable
function $f\ge 0$ on $\R^d$  with $f/h\le c$ for some $c>0$ and that
the set of  discontinuous points of $f$ has zero $m$-measure, we have
\begin{eqnarray}\label{limit-Xh3}
\lim_{t\rightarrow\infty}e^{-\lambda^h_0 t}\langle f/h, X^h_t\rangle (\omega)& =&
W_\infty(\omega)\int_{\R^d}(f/h)(y)m^h(dy)\nonumber\\
&=&W_\infty(\omega)\Big(\frac{c}{\pi}
\Big)^{d/2}\int_{\R^d}f(y)e^{(\upsilon-c)|y|^2}\,dy,
\end{eqnarray}
where $W_\infty(\omega)$ is the limit of the martingale $W_t:=e^{-\lambda_ct}
\langle 1, X_t^h\rangle=e^{-\lambda_ct}\langle h, X_t\rangle$ as $t\to\infty$.
We rewrite \eqref{limit-Xh3} to get the limit result on $X$:
 \begin{eqnarray}\label{limit-X}
\lim_{t\rightarrow\infty}e^{-\lambda_c t}\langle f, X_t\rangle (\omega) &=&
W_\infty(\omega)\int_{\R^d}\left(\frac{c}{\pi}\right)^{d/2}e^{(\upsilon-c)
\|y\|^{2}}f(y)dy\nonumber\\
&=&W_\infty(\omega)
\int_{\R^d}\widetilde\phi_0(y)f(y)dy,
\end{eqnarray}where $\widetilde \phi_0=\left(\frac{c}{\pi}\right)^{d/2}
e^{(\upsilon-c)\|y\|^{2}}$.
Since  $h$ is bounded from below,  in the weak topology, $e^{-\lambda_c t}
X_t\to W_\infty(\omega)\widetilde \phi_0(x)dx$, $\mathbb P_\mu$-a.s.,
for any $\mu\in \mathcal M_F(\R^d)$. This example covers \cite[Example 4.7]{EKW14}.

}
\end{example}

\begin{singlespace}

\end{singlespace}

\vskip 0.2truein
\vskip 0.2truein

\noindent{\bf Zhen-Qing Chen:} Department of Mathematics,
University of Washington,
Seattle, WA 98195, USA.
Email: {\texttt zqchen@uw.edu}

\medskip

\noindent{\bf Yan-Xia Ren:} LMAM School of Mathematical Sciences \& Center for
Statistical Science, Peking
University,  Beijing, 100871, P.R. China. Email: {\texttt
yxren@math.pku.edu.cn}

\medskip

\noindent {\bf Renming Song:} Department of Mathematics,
University of Illinois,
Urbana, IL 61801, U.S.A.
Email: {\texttt rsong@math.uiuc.edu}

\medskip

\noindent{\bf Rui Zhang:} LMAM School of Mathematical Sciences, Peking
University,  Beijing, 100871, P.R. China. Email: {\texttt
ruizhang8197@gmail.com}

\end{doublespace}
\end{document}